\documentclass[journal]{IEEEtran}
\ifCLASSINFOpdf
\else
\fi
\usepackage{multirow}
\usepackage{amsmath}
\usepackage{amssymb}
\usepackage{url}
\usepackage{graphicx}
\usepackage{breakurl}
\usepackage[numbers,sort&compress]{natbib}
\begin{document}
\title{A Hybrid MILP and IPM for Dynamic Economic Dispatch with Valve Point Effect}

\author{Shanshan~Pan,
        Jinbao~Jian~
        and~Linfeng~Yang% <-this % stops a space
\thanks{This work was supported by the Natural Science Foundation of China (No.51407037); the Natural Science Foundation of Guangxi (No.2016GXNSFDA380019, No.2014GXNSFFA118001); the Open Project Program of Guangxi Colleges and Universities Key Laboratory of Complex System Optimization and Big Data Processing (No.2016CSOBDP0205).}% <-this % stops a space

\thanks{S.S. Pan, is with the College of Electrical Engineering, Guangxi University, Nanning 530004, China. (e-mail: sspan@mail.gxu.cn).}% <-this % stops a space

\thanks{J.B. Jian, the corresponding author, is with the Guangxi Colleges and Universities Key Laboratory of Complex System Optimization and Big Data Processing, Yulin Normal University, Yulin 537000, China(e-mail: jianjb@gxu.edu.cn).}% <-this % stops a space

\thanks{L.F. Yang is with the College of Computer Electronics and Information, Guangxi
University, Nanning 530004, China.(e-mail: ylf@gxu.edu.cn).}
}

\maketitle

\begin{abstract}
Dynamic economic dispatch with valve-point effect (DED-VPE) is a non-convex and non-differentiable optimization problem which is difficult to solve efficiently. In this paper, a hybrid mixed integer linear programming (MILP) and interior point method (IPM), denoted by MILP-IPM, is proposed to solve such a DED-VPE problem, where the complicated transmission loss is also included. Due to the non-differentiable characteristic of DED-VPE, the classical derivative-based optimization methods can not be used any more.
With the help of model reformulation, a differentiable non-linear programming (NLP) formulation which can be directly solved by IPM is derived. However, if the DED-VPE is solved by IPM in a single step, the optimization will easily trap in a poor local optima due to its non-convex and multiple local minima characteristics. To exploit a better solution, an MILP method is required to solve the DED-VPE without transmission loss, yielding a good initial point for IPM to improve the quality of the solution. Simulation results demonstrate the validity and effectiveness of the proposed MILP-IPM in solving DED-VPE.
\end{abstract}


% \fi
\IEEEpeerreviewmaketitle

%\section{Introduction}
%\IEEEPARstart{T}{his} demo file is intended to serve as a ``starter file''
%for IEEE journal papers produced under \LaTeX\ using
%IEEEtran.cls version 1.8b and later.
%% You must have at least 2 lines in the paragraph with the drop letter
%% (should never be an issue)
%I wish you the best of success.

\section{Introduction}

Dynamic economic dispatch (DED) is one of the fundamental issues for the optimal economic operation in power system, which aims at allocating the customers' load demands among
the available thermal power generating units in an economic, secure and reliable way at a certain time of interest \cite{Xia2010}. Traditionally, the generation cost function
of DED is approximated by a convex, quadratic and differentiable polynomials. However, in actual operation, wire drawing effects, occurring as each steam admission valve in a
turbine starts to open, produce a rippling effect on the generation cost curve\cite{GA1993}, which is known as the valve-point effect (VPE). To model the effect of
valve-points, a recurring rectified sinusoid contribution is added to the input-output equation \cite{GA1993}, which makes the generation cost function non-convex,
non-differentiable and multiple extremal. When VPE is ignored, the rough approximation of the generation cost function will introduce some inaccuracies into the dispatch results. In order to improve the optimality of the solution, a more accurate DED model, dynamic economic dispatch with
valve-point effect (DED-VPE) should be considered. However, when VPE is considered, some non-convex, non-differentiable and multiple extremal characteristics are introduced, which makes the solution for DED-VPE more challenging.

In the past decades, a number of optimization methods have been proposed to solve the DED-VPE. Due to the intractability of the problem, most of the currently available approaches for DED-VPE are heuristic optimization
techniques\cite{GA1993,EPSQP2002,DGPSO2005,SA2006,APSO2008,DE2008,IPSO2009,ICPSO2010,AHDE2010,CSAPSO2011,ABC2011,AIS2011,ECE2011,HHS2011,CDE2011,CSDE2011,EAPSO2012,
TVACIPSO2012,ICA2012,EBSO2013,MTLA2013,CDBCO2014,TLA2015,CSO2015,EA2016},
such as genetic algorithm (GA) \cite{GA1993}, evolutionary programming (EP) \cite{EPSQP2002}, simulated annealing (SA) \cite{SA2006}, particle swarm optimization (PSO)
\cite{APSO2008}, differential evolution algorithm (DE) \cite{DE2008}, artificial bee colony algorithm (ABC) \cite{ABC2011}, artificial immune system (AIS) \cite{AIS2011},
enhanced cross-entropy (ECE) \cite{ECE2011}, harmony search (HS) \cite{HHS2011}, imperialist competitive algorithm (ICA) \cite{ICA2012}, bee swarm optimization algorithm (BSO) \cite{EBSO2013}, teaching-learning algorithm (TLA) \cite{TLA2015}, etc.
These heuristic optimization techniques are population-based search methods which do not depend on the gradient and Hessian operators of the objective function. So, they can be applied to solve the DED-VPE problem effectively.
However, they are quite sensitive to various parameter settings and solution may be different at each trial. Hence, hybrid methods which combine several heuristic techniques or deterministic approaches \cite{EPSQP2002,MHEPSQP2005,AISSQP2009,SOASQP2010,HIGA2013,HBPSO2014} such as hybrid evolutionary programming and sequential quadratic
programming (EP-SQP) \cite{EPSQP2002}, hybridization of artificial immune systems and sequential quadratic programming (AIS-SQP) \cite{AISSQP2009}, hybrid seeker optimization algorithm and sequential quadratic programming (SOA-SQP) \cite{SOASQP2010}, hybrid immune-genetic algorithm (HIGA) \cite{HIGA2013}, etc, tend to be more efficient than the individual methods. However, they still have the intrinsic drawbacks of the heuristic method we mentioned above.

Unlike heuristics, deterministic mathematical programming-based optimization techniques can solve to a robust result due to the solid mathematical foundations and the
availability of the powerful software tools. Therefore,
a strategy recently appeared for DED-VPE is to reformulate the generation cost function, yielding a good optimization model that can be solved by a deterministic method. In
\cite{MIQP2014}, the non-convex and non-differentiable cost function caused by VPE is piecewise linearized, then a mixed integer quadratic programming (MIQP) method can be used to
solve DED-VPE. But when the MIQP formulation is directly solved by using a mixed integer programming (MIP) solver, the optimization will suffer convergence stagnancy and run out of memory. As a result, the multi-step method, the warm start technique and the range restriction scheme are required. However, the range restriction scheme just restricts the solution space to a subspace where the global optimal solution would probably lie in. Consequently, the optimality of the solution for the MIQP can not be guaranteed. Besides, when the complicated transmission loss is considered, base on the above process, more tedious adjustment techniques are needed. Different from \cite{MIQP2014}, the whole generation cost function is considered for piecewise linearization in \cite{MILP2017}, then a mixed integer linear programming (MILP) formulation is proposed to solve the DED-VPE, where a global optimal solution within a preset tolerance can be guaranteed. But the transmission loss is not considered in \cite{MILP2017}. Therefore, more efforts are worth developing effective deterministic mathematical programming-based optimization methods for DED-VPE to obtain better dispatch results.

In this paper, a hybrid deterministic method that integrates the MILP and interior point method (IPM), denoted by MILP-IPM, is proposed for solving the DED-VPE problem while transmission loss is included. Due to the non-differentiable characteristic of DED-VPE, the classical derivative-based optimization methods can not be used any more. With the help of model reformulation, we derive a non-linear programming (NLP) formulation of DED-VPE, which can be solved by the polynomial time IPM immediately. However, IPM is a local optimization method. If the DED-VPE is solved by IPM in a single step, the optimization will easily trap in a poor local optima due to its non-convex and multiple local minima characteristics. In order to overcome this deficiency, MILP method \cite{MILP2017} is combined to generate a good initial point for IPM. And then, solving its NLP formulation via IPM, a good optimal solution for DED-VPE can thus be determined.

The rest of this paper is organized as follows. Section 2 describes the mathematical formulation of DED-VPE. Section 3 derives an NLP formulation for DED-VPE. Section 4 introduces the implementation of MILP-IPM. Section 5 presents the simulation results and analysis. Section 6 draws the conclusions.

\section{Mathematical formulation of DED-VPE}
\label{DED-VPE}

The DED-VPE problem usually can be formulated as a non-convex and non-differentiable optimization problem. The objective of DED-VPE is to minimize the total generation cost over a scheduled time horizon, which can be written as:
\begin{equation}\label{objective}
\min~\sum\limits_{t=1}^T\sum\limits_{i=1}^N c(P_{i,t})
\end{equation}
where
\begin{equation}\label{cost}
c(P_{i,t})=\alpha_i+\beta_iP_{i,t}+\gamma_iP_{i,t}^2+e_i|\sin(f_i(P_{i,t}-P_i^{\min}))|
\end{equation}
where $P_{i,t}$ is the power output of uint $i$ in period $t$; $P_i^{\min}$ is the minimum power output of unit $i$; $\alpha_i$, $\beta_i$, $\gamma_i$, $e_i$ and $f_i$ are positive coefficients of unit $i$; $N$ and $T$ are total number of units and periods, respectively.

The minimized DED-VPE problem should be subjected to the constraints as follows.

\begin{itemize}
\item Power balance equations
\begin{equation}\label{balance}
\sum\limits_{i=1}^N P_{i,t}=D_{t}+P_{t}^{loss}, ~~~\forall~t
\end{equation}
where $D_{t}$ is the load demand in period $t$; $P_{t}^{loss}$ is the transmission loss in period $t$, which can be calculated based on the B-Matrix loss coefficients and expressed in the quadratic form as given below \cite{Saadat1999}
\begin{equation}\label{loss}
P_{t}^{loss}=\sum\limits_{i=1}^N\sum\limits_{j=1}^N P_{i,t}B_{i,j}P_{j,t}, ~~~\forall~t
\end{equation}
where $B_{i,j}$ is the $(i,j)$-th element of the transmission loss coefficient matrix.
\item Power generation limits
\begin{equation}\label{G-limits}
P_i^{\min}\le{P_{i,t}}\le{P_i^{\max}}, ~~~\forall~i,t
\end{equation}
where $P_i^{\max}$ is the maximum power output of unit $i$.
\item Ramp rate limits
\begin{equation}\label{R-limits}
DR_{i}\le{P_{i,t}-P_{i,t-1}}\le{UR_{i}}, ~~~\forall~i,t
\end{equation}
where $DR_{i}$ and $UR_{i}$ are the ramp-down and ramp-up rates of unit $i$, respectively.
\item Spinning reserve constraints
\begin{equation}\label{SR1}
\left\{
\begin{aligned}
&SR_{i,t}\le P_{i}^{\max}-P_{i,t},     \\
&SR_{i,t}\le\tau UR_{i},~~\forall~i,t \\
\end{aligned}
\right.
\end{equation}
\begin{equation}\label{SR2}
\sum\limits_{i=1}^N SR_{i,t}\geq R_{t}, ~~~\forall~t
\end{equation}
where $SR_{i,t}$ is the spinning reserve provided by unit $i$ in period $t$; $R_{t}$ is the system spinning reserve requirement in period $t$; $\tau$ is the time duration for units to deliver reserve \cite{MTLA2013,MIQP2014}.
\end{itemize}

\section{An NLP formulation for DED-VPE}

As we can see from the section \ref{DED-VPE}, DED-VPE is a non-convex and non-differentiable optimization problem which is hard to tackle. Due to the non-differentiable characteristic of DED-VPE, the classical mathematical programming-based methods, also known as derivative-based optimization
methods, are not suitable any more. To overcome this difficulty, we replace $|\sin(f_i(P_{i,t}-P_i^{\min}))|$ for (\ref{cost}) with an auxiliary variable $s_{i,t}$, then the objective function (\ref{objective}) can be equivalent to
\begin{eqnarray}
\min&\sum\limits_{t=1}^T\sum\limits_{i=1}^N (\alpha_i+\beta_iP_{i,t}+\gamma_iP_{i,t}^2+e_is_{i,t})   \\
\label{auxiliary1}
s.t.&s_{i,t} \geq \sin(f_i(P_{i,t}-P_i^{\min}))  ~~~~~~~~~   \\
\label{auxiliary2}
&s_{i,t} \geq -\sin(f_i(P_{i,t}-P_i^{\min})) ~~~~~~~
\end{eqnarray}
By introducing some slack variables $u_{i,t}^0\geq 0$, $v_{i,t}^0\geq 0$, the new inequality constraints (\ref{auxiliary1}) and (\ref{auxiliary2}) are converted into equality constraints
\begin{eqnarray}
\label{slack1}
s_{i,t}-\sin(f_i(P_{i,t}-P_i^{\min}))- u_{i,t}^0=0     \\
\label{slack2}
s_{i,t}+\sin(f_i(P_{i,t}-P_i^{\min}))- v_{i,t}^0=0
\end{eqnarray}
From (\ref{slack1}) and (\ref{slack2}), we can get the following equations
\begin{eqnarray}
\label{slack3}
2s_{i,t}-u_{i,t}^0-v_{i,t}^0=0  \\
\label{slack4}
-2\sin(f_i(P_{i,t}-P_i^{\min}))- u_{i,t}^0+v_{i,t}^0=0
\end{eqnarray}
Then we have
\begin{eqnarray}
\label{slack5}
s_{i,t}-\frac{u_{i,t}^0}{2}-\frac{v_{i,t}^0}{2}=0  \\
\label{slack6}
\sin(f_i(P_{i,t}-P_i^{\min}))+\frac{u_{i,t}^0}{2}-\frac{v_{i,t}^0}{2}=0
\end{eqnarray}
Let $u_{i,t}=\frac{u_{i,t}^0}{2}$, $v_{i,t}=\frac{v_{i,t}^0}{2}$, we can obtain that
\begin{eqnarray}
\label{slack01}
s_{i,t}-u_{i,t}-v_{i,t}=0   \\
\label{slack02}
\sin(f_i(P_{i,t}-P_i^{\min}))+u_{i,t}-v_{i,t}=0 \\
\label{slack03}
u_{i,t}\geq 0,~~v_{i,t}\geq 0
\end{eqnarray}

Consequently, the DED-VPE can be formulated as the following differentiable NLP formulation.

\begin{equation}\label{NLP}
\begin{aligned}
\min~&\sum\limits_{t=1}^T\sum\limits_{i=1}^N(\alpha_i+\beta_i P_{i,t}+\gamma_i P_{i,t}^2+e_is_{i,t})\\
s.t.~&(\ref{balance})-(\ref{SR2}), (\ref{slack01})-(\ref{slack03})
\end{aligned}
\end{equation}

\section{Implementation of MILP-IPM}

It is well known that, IPM is a powerful tool for solving non-linear optimization problem and it has successfully been implemented
in various issues in power system \cite{Wu2001}, such as optimal power flows, state estimation, hydrothermal coordination, economic dispatch, etc. Hence, the above NLP formulation
(\ref{NLP}) for DED-VPE can be directly solved by using IPM. However, if the NLP formulation of DED-VPE is solved by IPM in a single step, the solution will easily trap in a poor local optima due to its non-convex, non-linear and multiple local minima characteristics.

To exploit a better solution, an MILP method presented in \cite{MILP2017} is incorporated to find a good initial point for IPM.
In \cite{MILP2017}, when transmission loss is not considered, the DED-VPE is reformulated into an MILP formulation which can be solved by a state-of-the-art MIP solver directly and efficiently. Hence, a global optimal solution within a preset tolerance, which can be used as an initial point for IPM to improve the quality of the eventual dispatch result, can be guaranteed via an enumeration algorithm.

\subsection{An MILP formulation for DED-VPE}

To obtain an MILP formulation of DED-VPE, ${L_i} + 1$ break points are chosen over a generation interval $[P_i^{\min },P_i^{\max }]$, such that $ P_i^{min}=a_{0,i}\le a_{1,i}\le \cdots \le a_{{L_i},i}=P_i^{max}$. Segment variables $P_{l,i,t}$ and binary variables $z_{l,i,t}$ $(l =1,\cdots {L_i})$ are introduced to make $P_{m,i,t}=P_{i,t}$ and $P_{l,i,t}=0(l\ne m)$ when the $P_{i,t}$ lies in segment $m(m\in \{1,\cdots {L_i}\} )$. The detail of the MILP formulation for DED-VPE without considering the transmission loss is expressed as \cite{MILP2017}:

\begin{equation}\label{MILP}
\begin{aligned}
\min~&\sum\limits_{t=1}^T\sum\limits_{i=1}^N\sum\limits_{l=1}^{L_i}(k_{l,i}P_{l,i,t}+b_{l,i}z_{l,i,t})\\
s.t.~&P_{i,t}=\sum\limits_{l = 1}^{L_i}{P_{l,i,t}}  \\
&a_{l-1,i}z_{l,i,t}\le{P_{l,i,t}}\le{a_{l,i}}z_{l,i,t}  \\
&\sum\limits_{l=1}^{L_i}{z_{l,i,t}}=1\\
&z_{l,i,t}\in\{0,1\} \\
&(\ref{balance}), (\ref{G-limits})-(\ref{SR2})
\end{aligned}
\end{equation}
where $P_{t}^{loss}=0$ for (\ref{balance}) and $L_i$, $k_{l,i}$, $b_{l,i}$ are calculated as follows
\begin{equation}\label{lkb1}
\left\{
\begin{aligned}
&L_i=ceil(M\frac{f_i(P_i^{max}-P_i^{min})}\pi)  \\
&k_{l,i}=\frac{c(a_{l,i})-c(a_{l-1,i})}{a_{l,i}-a_{l-1,i}}  \\
&b_{l,i}=c(a_{l - 1,i})-k_{l,i}a_{l - 1,i}.  \\
\end{aligned}
\right.
\end{equation}
Above $ceil(x)$ means round $x$ to the nearest integer greater than or equal to $x$ and $M$ is the number of equal segments on each $sin(x)$ where $x$ belongs to $[0,\pi]$.

\subsection{Details of MILP-IPM}

Now the MILP-IPM for DED-VPE as a whole is summarized as follows.

\textbf{Step 1}: Solve the MILP formulation (\ref{MILP}) by using MILP method to obtain a global optimal solution within a preset tolerance for DED-VPE without transmission loss.

\textbf{Step 2}: Solve the NLP formulation (\ref{NLP}) by using IPM, where the initial point is taken at the solution gained in Step 1, to obtain a good local optimal solution for DED-VPE with transmission loss.

In step 1, an MILP formulation for DED-VPE without transmission loss is solved, yielding a solution which is close to a good optimal solution for DED-VPE with transmission loss. This is mainly because, in a DED problem, the transmission loss at each period is small compared to the corresponding load demand \cite{MIQP2014}. When transmission loss is ignored, the solution gained in step 1 is ``the most" economic. When transmission loss is considered, it means that more outputs are needed and some transmission loss constraints will be added into the original model. Then base on the most economic initial solution, some units outputs will be fine tuned (since the transmission loss is small compared to the load demand ) by IPM in step 2, to meet the new constraints and attain a new economic solution. The following simulation results conform well with the above analysis.

\section{Simulation results and analysis}

To assess the validity and effectiveness of the proposed MILP-IPM, in this section, several test systems which are widely studied for DED-VPE with and without consideration of transmission loss constraints over a scheduled time horizon of 24 h are simulated. For a fairer comparison with most of the existing methods, the spinning reserve constraints (\ref{SR1}) and (\ref{SR2}) are not included in our simulations. Actually, it does not reduce the difficulty of the problem since the spinning reserve constraints are some simple linear constraints which is very easy to handle for MILP and IPM.

Since the computation time highly depends on the computer system used, so the CPU execution times for different methods may not be directly comparable due to different computers used. In order to have a fair comparison regarding the computational effort, the CPU chip frequency from the used computer is used to convert the CPU times obtained from different methods into a common base for comparison \cite{scaletime2008}:
\begin{equation}\label{Stime}
\begin{aligned}
&Scaled~CPU~time\\
&=\frac{Given~CPU~speed}{Base~CPU~speed} Given~CPU~time.
\end{aligned}
\end{equation}
In this paper, the base CPU speed is 2.4 GHz and we denote $Scaled~CPU~time$ by S-time for short.

In this section, we first carry out two experiments on DED-VPE without transmission loss to compare MILP-IPM with MILP and other methods, showing the potential of MILP-IPM for seeking better solutions. Then two experiments on DED-VPE with transmission loss are carried out to demonstrate the validity and  effectiveness of MILP-IPM in solving DED-VPE. All cases are performed on an Intel Core 2.5 GHz Dell-notebook with 8 GB of RAM. Our models are coded in YALMIP \cite{Lofberg2004} within the MATLAB R2014a and optimised using CPLEX 12.6.1 \cite{CPLEX2015} for solving the MILP (\ref{MILP}) and IPOPT 3.12.6 \cite{IP2006} for solving the NLP (\ref{NLP}).

\subsection{DED-VPE without transmission loss}

\textbf{Case 1}: A 5-unit system without transmission loss

The first test case is incorporated five thermal units. The characteristics of the thermal units and load demands are taken from \cite{SA2006}. Firstly, we directly solve its MILP formulation using CPLEX to 3.2\% optimality. In this formulation, segment parameter $M$ is set to 4. After the optimization, an optimal solution is yielded in 0.82 min with a total generation cost $42563\$$. Taking such a solution as an initial point in step 2,  we solve its NLP formulation using IPOPT with the default options, then a local optimal solution with an optimal value $42524\$$ is found in 0.04 min.

Table \ref{cost5} lists comparison results of the total generation cost obtained by MILP and MILP-IPM. It is clearly seen that the proposed MILP-IPM can make a certain improvement for 5-unit system.

\begin{table}[!ht]
\caption{\label{cost5} Summary results for the 5-unit system without loss}
\centering
\resizebox{0.46\textwidth}{!}{
\begin{tabular}{ccccc}
\hline
\multirow{2}{*}{Method} & \multicolumn{3}{c}{Total generation cost (\$)}&\multirow{2}{*}{S-time(min)}\\
\cline{2-4}
& Minimum & Average & Maximum & \\
\hline
MILP     &    42563 &  &&  0.82\\
MILP-IPM &    42524 &  &&  0.86\\
\hline
\end{tabular}}
\end{table}

The outputs obtained by using MILP-IPM for the 5-unit system are given in Table \ref{output5} for verification.

\begin{table}[!ht]
\caption{\label{output5} Outputs (MW) for the 5-unit system without loss}
\centering
\resizebox{0.42\textwidth}{!}{
\begin{tabular}{c|ccccc}
\hline
$t$  & unit $1$&  unit $2$&  unit $3$&  unit $4$&  unit $5$\\
\hline
 1& 16.7925& 98.5398& 30.0000& 124.9079& 139.7598\\
\hline
 2& 10.0000& 98.5398& 61.7925& 124.9079& 139.7598\\
\hline
 3& 10.0000& 98.5398& 101.7925& 124.9079& 139.7598\\
\hline
 4& 10.0000& 98.5398& 112.6735& 124.9079& 183.8788\\
\hline
 5& 10.0000& 80.8990& 112.6735& 124.9079& 229.5196\\
\hline
 6& 40.0000& 100.5398& 112.6735& 125.2671& 229.5196\\
\hline
 7& 10.0000& 98.5398& 112.6735& 175.2671& 229.5196\\
\hline
 8& 10.0000& 91.9911& 112.6735& 209.8158& 229.5196\\
\hline
 9& 39.4513& 98.5398& 112.6735& 209.8158& 229.5196\\
\hline
10& 53.4513& 98.5398& 112.6735& 209.8158& 229.5196\\
\hline
11& 69.4513& 98.5398& 112.6735& 209.8158& 229.5196\\
\hline
12& 75.0000& 112.9911& 112.6735& 209.8158& 229.5196\\
\hline
13& 53.4513& 98.5398& 112.6735& 209.8158& 229.5196\\
\hline
14& 39.4513& 98.5398& 112.6735& 209.8158& 229.5196\\
\hline
15& 34.0000& 98.5398& 112.6735& 179.2671& 229.5196\\
\hline
16& 10.0000& 98.5398& 112.6735& 129.2671& 229.5196\\
\hline
17& 10.0000& 80.8990& 112.6735& 124.9079& 229.5196\\
\hline
18& 10.0000& 98.5398& 112.6735& 157.2671& 229.5196\\
\hline
19& 10.0000& 94.5398& 112.6735& 207.2671& 229.5196\\
\hline
20& 40.0000& 111.9911& 112.6735& 209.8158& 229.5196\\
\hline
21& 29.4513& 98.5398& 112.6735& 209.8158& 229.5196\\
\hline
22& 10.0000& 96.7800& 108.8846& 209.8158& 179.5196\\
\hline
23& 10.0000& 98.5398& 68.8846& 209.8158& 139.7598\\
\hline
24& 10.0000& 73.4244& 30.0000& 209.8158& 139.7598\\
\hline
\end{tabular}}
\end{table}

\begin{table}[!ht]
\caption{\label{cost10} Summary results for the 10-unit system without loss}
\centering
\resizebox{0.46\textwidth}{!}{
\begin{tabular}{ccccc}
\hline
\multirow{2}{*}{Method} & \multicolumn{3}{c}{Total generation cost (\$)}& \multirow{2}{*}{S-time(min)}\\
\cline{2-4}
& Minimum & Average & Maximum & \\
\hline
SQP \cite{EPSQP2002}        &1051163  &NA       &NA        &0.42 \\
EP \cite{EPSQP2002}         &1048638  &NA       &NA        &15.05 \\
CDE \cite{CDE2011}          &1036756  &1040586  &1452558   &0.20\\
GA \cite{ABC2011}           &1033481  &1038014  &1042606   &3.59 \\
EP-SQP \cite{EPSQP2002}     &1031746  &1035748  &NA        &7.26 \\
AIS-SQP\cite{AISSQP2009}    &1029900  &NA       &NA        &NA\\
MHEP-SQP\cite{MHEPSQP2005}  &1028924  &1031179  &NA        &21.23\\
DGPSO \cite{DGPSO2005}      &1028835  &1030183  &NA        &4.81\\
PSO \cite{ABC2011}          &1027679  &1031716  &1034340   &3.85 \\
SOA \cite{SOASQP2010}       &1023946  &1026289  &1029213   &NA  \\
IPSO \cite{IPSO2009}        &1023807  &1026863  &NA        &0.05\\
CSDE \cite{CSDE2011}        &1023432  &1026475  &1027634   &0.3\\
CE \cite{ECE2011}           &1022702  &1024024  &NA        &0.33 \\
ECE  \cite{ECE2011}         &1022272  &1023334  &NA        &0.33\\
AIS \cite{AIS2011}          &1021980  &1023156  &1024973   &25.35 \\
ABC\cite{ABC2011}           &1021576  &1022686  &1024316   &3.47 \\
CDBCO \cite{CDBCO2014}      &1021500  &1024300  &NA        &0.73\\
SOA-SQP\cite{SOASQP2010}    &1021460  &1023841  &1026852   &NA\\
AHDE \cite{AHDE2010}        &1020082  &1022474  &1024484   &1.10\\
CDE\cite{CDE2011}           &1019123  &1020870  &1023115   &0.32\\
HHS  \cite{HHS2011}         &1019091  &NA       &NA        &10.19\\
ICPSO \cite{ICPSO2010}      &1019072  &1020027  &NA        &0.35\\
CSAPSO\cite{CSAPSO2011}     &1018767  &1019874  &NA        &0.350\\
EAPSO \cite{EAPSO2012}      &1018510  &1018701  &1019302   &0.63\\
HIGA \cite{HIGA2013}        &1018473  &1019328  &1022284   &4.41\\
ICA \cite{ICA2012}          &1018467  &1019291  &1021796   &NA\\
TVAC-IPSO\cite{TVACIPSO2012}&1018217  &1018965  &1020418   &2.72\\
HBPSO\cite{HBPSO2014}       &1018159  &1019850  &1021813   &3.09\\
CSO \cite{CSO2015}          &1017660  &1018120  &1019286   &0.90\\
EBSO \cite{EBSO2013}        &1017147  &1017526  &1017891   &0.15\\
MILP     &1016316  & && 0.94   \\
MILP-IPM &1016311  & && 1.02 \\
\hline
\multicolumn{5}{c}{\scriptsize{NA denotes that the value was not available in the literature.}}
\end{tabular}}
\end{table}

\begin{table*}[!ht]
\caption{\label{output10} Outputs (MW) for the 10-unit system without loss}
\centering
\resizebox{0.8\textwidth}{!}{
\begin{tabular}{c|cccccccccc}
\hline
$t$  & unit $1$&  unit $2$&  unit $3$&  unit $4$&  unit $5$& unit $6$&  unit $7$&  unit $8$&  unit $9$&  unit $10$\\
\hline
1& 150.0000& 142.2665& 186.8267& 60.0000& 122.8666& 122.4498& 129.5904& 47.0000& 20.0000& 55.0000\\
\hline
 2& 150.0000& 222.2665& 230.6932& 60.0000& 73.0000& 122.4498& 129.5904& 47.0000& 20.0000& 55.0000\\
\hline
 3& 150.0000& 302.2665& 298.6932& 60.0000& 73.0000& 122.4498& 129.5904& 47.0000& 20.0000& 55.0000\\
\hline
 4& 226.6242& 316.7994& 305.6696& 60.0000& 122.8666& 122.4498& 129.5904& 47.0000& 20.0000& 55.0000\\
\hline
 5& 226.6242& 396.7994& 299.6696& 60.0000& 122.8666& 122.4498& 129.5904& 47.0000& 20.0000& 55.0000\\
\hline
 6& 303.2484& 396.7994& 297.3995& 60.0000& 172.7331& 146.2292& 129.5904& 47.0000& 20.0000& 55.0000\\
\hline
 7& 379.8726& 396.7994& 303.3995& 66.8430& 172.7331& 122.4498& 129.5904& 55.3121& 20.0000& 55.0000\\
\hline
 8& 379.8726& 396.7994& 297.3995& 116.8430& 172.7331& 122.4498& 129.5904& 85.3121& 20.0000& 55.0000\\
\hline
 9& 456.4968& 396.7994& 297.3995& 131.6073& 222.5997& 122.4498& 129.5904& 90.0000& 22.0571& 55.0000\\
\hline
10& 456.4968& 396.7994& 297.8493& 181.6073& 222.5997& 160.0000& 129.5904& 120.0000& 52.0571& 55.0000\\
\hline
11& 456.4968& 396.7994& 321.8493& 231.6073& 222.5997& 160.0000& 129.5904& 120.0000& 52.0571& 55.0000\\
\hline
12& 456.4968& 460.0000& 332.6487& 231.6073& 222.5997& 160.0000& 129.5904& 120.0000& 52.0571& 55.0000\\
\hline
13& 456.4968& 396.7994& 297.8493& 181.6073& 222.5997& 160.0000& 129.5904& 120.0000& 52.0571& 55.0000\\
\hline
14& 456.4968& 396.7994& 297.3995& 131.6073& 222.5997& 122.4498& 129.5904& 90.0000& 22.0571& 55.0000\\
\hline
15& 379.8726& 396.7994& 312.5546& 110.0000& 172.7331& 122.4498& 129.5904& 77.0000& 20.0000& 55.0000\\
\hline
16& 303.2484& 396.7994& 297.0454& 60.0000& 122.8665& 122.4498& 129.5904& 47.0000& 20.0000& 55.0000\\
\hline
17& 226.6242& 396.7994& 299.6696& 60.0000& 122.8666& 122.4498& 129.5904& 47.0000& 20.0000& 55.0000\\
\hline
18& 303.2484& 396.7994& 306.0237& 66.8430& 172.7331& 122.4498& 129.5904& 55.3121& 20.0000& 55.0000\\
\hline
19& 379.8726& 396.7994& 297.3995& 116.8430& 172.7331& 122.4498& 129.5904& 85.3121& 20.0000& 55.0000\\
\hline
20& 456.4968& 460.0000& 332.5857& 120.4152& 222.5997& 160.0000& 129.5904& 85.3121& 50.0000& 55.0000\\
\hline
21& 456.4968& 389.5329& 297.3995& 119.5092& 222.5997& 148.5593& 129.5904& 85.3121& 20.0000& 55.0000\\
\hline
22& 379.8726& 309.5329& 283.9997& 69.5092& 172.7331& 122.4498& 129.5904& 85.3121& 20.0000& 55.0000\\
\hline
23& 303.2484& 229.5329& 203.9997& 60.0000& 122.8665& 122.4498& 129.5904& 85.3121& 20.0000& 55.0000\\
\hline
24& 226.6242& 222.2665& 189.7569& 60.0000& 73.0000& 122.4498& 129.5904& 85.3121& 20.0000& 55.0000\\
\hline
\end{tabular}}
\end{table*}

\begin{table}[!ht]
\caption{\label{cost5withloss} Summary results for the 5-unit system with loss}
\centering
\resizebox{0.46\textwidth}{!}{
\begin{tabular}{ccccc}
\hline
\multirow{2}{*}{Method} & \multicolumn{3}{c}{Total generation cost (\$)}& \multirow{2}{*}{S-time(min)} \\
\cline{2-4}
& Minimum & Average & Maximum & \\
\hline
SA\cite{SA2006}             &  47356&     NA&       NA&    4.40\\
GA \cite{ABC2011}           &  44862&  44922&    45894&    4.43\\
APSO\cite{APSO2008}         &  44678&     NA&	    NA&    NA\\
AIS \cite{AIS2011}          &  44385&  44759&    45554&    5.33\\
PSO \cite{ABC2011}          &  44253&  45657&    46403&    4.73\\
ABC \cite{ABC2011}          &  44046&  44065&    44219&    4.39\\
EAPSO\cite{EAPSO2012}       &  43784&     NA&       NA&    0.42\\
HBPSO\cite{HBPSO2014}       &  43223&  43732&	 44252&	   1.54\\
DE \cite{DE2008}            &  43213&  43813&    44247&    6.00 \\
HHS \cite{HHS2011}          &  43155&     NA&       NA&    2.33\\
TVAC-IPSO\cite{TVACIPSO2012}&  43137&  43186&    43302&    1.07\\
HIGA \cite{HIGA2013}        &  43125&  43162&    43259&    2.06\\
ICA \cite{ICA2012}          &  43117&  43144&    43210&    NA \\
IPM     & 43443       &    &&  0.05\\
MILP-IPM & 43084      &    &&  0.87\\
\hline
\end{tabular}}
\end{table}

\textbf{Case 2}: A 10-unit system without transmission loss

The second test case contains ten thermal units. The characteristics of the thermal units and load demands are taken from \cite{EPSQP2002}. Firstly, we directly solve its MILP formulation using CPLEX to 0.3\% optimality. After the optimization, an optimal solution is yielded in 0.94 min with a total generation cost $1016316\$$. Using such a solution as an initial point in step 2, we solve its NLP formulation using IPOPT with the default options, then a local optimal solution with optimal value $1016311\$$ is found in 0.08 min.

The comparison results of the total generation cost obtained by using MILP-IPM and other methods are shown in Table \ref{cost10}. It is obvious that the proposed MILP-IPM can solve to the lowest generation cost among all the methods in a reasonable time. But at the same time, we also see that, although MILP-IPM can exploit a better solution in comparison with MILP, but the improvement is small. This is mainly because the solution achieved in step 1 is very well.

The outputs obtained by using MILP-IPM for the 10-unit system are given in Table \ref{output10} for verification.

\subsection{DED-VPE with transmission loss}

Since transmission loss can not be avoided in a power distribution system and it is critical in real-world DED problem, solution for DED-VPE with transmission loss has more values. Due to the data unavailability \cite{EA2016}, in our simulations, only two cases (5- and 10- unit systems) for DED-VPE with transmission loss are considered.

\textbf{Case 3}: A 5-unit system with transmission loss

The third test case is a 5-unit system in which the characteristics of the thermal units and load demands are the same as those in case 1. In this case, the transmission loss is considered. Owing to the limits of space, the loss coefficients are not listed here. One can refer to \cite{SA2006}.

Since the transmission loss is included, unlike case 1, it can not be solved by MILP immediately. Fortunately, with the help of model reformulation, we derive a differentiable NLP formulation (\ref{NLP}) of DED-VPE, which can be solved by the powerful IPM directly.

In order to demonstrate the validity and effectiveness of the proposed NLP formulation and MILP-IPM for the DED-VPE with transmission loss, on the one hand, we directly solve its NLP formulation using IPOPT with the default options. Then a local optimal solution with an optimal value $43443\$$ is obtained in 0.05 min. On the other hand, we solve this system by using MILP-IPM, where the parameters in step 1 are set the same as the case 1. After the optimization, a local optimal solution with an optimal value $43084\$$ is found in 0.87 min. The results are compared with other methods in Table \ref{cost5withloss}.

As we can see in Table \ref{cost5withloss}, although the proposed NLP formulation can be solved to a local optimal solution in a short time, but when the DED-VPE is directly solved by IPM in a single step, the obtaining solution is not the best due to its non-convex and multiple local minima characteristics. While the MILP-IPM which combines MILP with IPM is employed, it can solve to the lowest generation cost among all the methods in a reasonable time.

But, we should note that, not only the optimality but also the feasibility should be considered for assessing the quality of the solution. In Table \ref{output5withloss}, the generation dispatch results obtained by using MILP-IPM for the 5-unit system are presented, where $\Delta P_{t}$ denotes the  violation degree of the power balance constraint at interval $t$, which is calculated by
\begin{equation}
\Delta P_{t}=|\sum\limits_{i=1}^N P_{i,t}-D_{t}-\sum\limits_{i=1}^N\sum\limits_{j=1}^N P_{i,t}B_{i,j}P_{j,t}|.   \nonumber
\end{equation}
Note that, the outputs which have been round-off are used for this calculation (i.e. the outputs in Table \ref{output5withloss}). Actually, when the original outputs are adopted, all the violations are less than 7e-7. In other words,
the solution obtained by using MILP-IPM strictly satisfies the power balance constraints at the same time. This is mainly the result of the rigorous theoretical foundations of interior point algorithm.

\begin{table}[!ht]
\caption{\label{output5withloss} Dispatch results (MW) for the 5-unit system}
\centering
\resizebox{0.48\textwidth}{!}{
\begin{tabular}{c|cccccccc}
\hline
$t$  & unit $1$&  unit $2$&  unit $3$&  unit $4$&  unit $5$&  Loss & $\Delta P_{t}$ \\
\hline
1& 20.6080& 98.5398& 30.0000& 124.9079& 139.7598& 3.8155 & 0.0000 \\
\hline
 2& 10.0000& 97.8835& 66.5747& 124.9079& 139.7598& 4.1259 & 0.0000 \\
\hline
 3& 10.0000& 98.5398& 106.5747& 124.9079& 139.7598& 4.7822 & 0.0000 \\
\hline
 4& 10.0960& 98.5398& 112.6735& 124.9079& 189.7598& 5.9770 & 0.0000 \\
\hline
 5& 10.0000& 87.5816& 112.6735& 124.9079& 229.5196& 6.6825 & 0.0001 \\
\hline
 6& 40.0000& 99.9508& 112.6735& 133.7152& 229.5196& 7.8590 & 0.0001 \\
\hline
 7& 10.0000& 98.5398& 112.6735& 183.7152& 229.5196& 8.4480 & 0.0001 \\
\hline
 8& 12.7090& 98.5398& 112.6735& 209.8158& 229.5196& 9.2577 & 0.0000 \\
\hline
 9& 42.7090& 105.4824& 112.6735& 209.8158& 229.5196& 10.2003 & 0.0000 \\
\hline
10& 64.0108& 98.5398& 112.6735& 209.8158& 229.5196& 10.5595 & 0.0000 \\
\hline
11& 75.0000& 104.0359& 112.6735& 209.8158& 229.5196& 11.0448 & 0.0000 \\
\hline
12& 75.0000& 124.7111& 112.6735& 209.8158& 229.5196& 11.7200 & 0.0000 \\
\hline
13& 64.0108& 98.5398& 112.6735& 209.8158& 229.5196& 10.5595 & 0.0000 \\
\hline
14& 49.6196& 98.5398& 112.6735& 209.8158& 229.5196& 10.1683 & 0.0000 \\
\hline
15& 19.6196& 98.5398& 112.6735& 202.8589& 229.5196& 9.2113 & 0.0001 \\
\hline
16& 10.0000& 82.1494& 112.6735& 152.8589& 229.5196& 7.2013 & 0.0001 \\
\hline
17& 10.0000& 87.5816& 112.6735& 124.9079& 229.5196& 6.6825 & 0.0001 \\
\hline
18& 10.0000& 98.5398& 112.6735& 165.2180& 229.5196& 7.9509 & 0.0000 \\
\hline
19& 12.7090& 98.5398& 112.6735& 209.8158& 229.5196& 9.2577 & 0.0000 \\
\hline
20& 42.7090& 119.9393& 112.6735& 209.8158& 229.5196& 10.6572 & 0.0000 \\
\hline
21& 39.3529& 98.5398& 112.6735& 209.8158& 229.5196& 9.9016 & 0.0000 \\
\hline
22& 10.0000& 98.5398& 112.6735& 209.8158& 181.8844& 7.9136 & 0.0001 \\
\hline
23& 12.4371& 98.5398& 72.6735& 209.8158& 139.7598& 6.2261 & 0.0001 \\
\hline
24& 10.0000& 75.8153& 32.6735& 209.8158& 139.7598& 5.0644 & 0.0000 \\
\hline
\end{tabular}}
\end{table}

\textbf{Case 4}: A 10-unit system with transmission loss

The fourth test case is a 10-unit system in which the characteristics of the thermal units and load demands are the same as those in case 2 and the loss coefficients are taken from \cite{DGPSO2005}.

Similar to the case 3, on the one hand, we directly solve its NLP formulation using IPOPT with the default options. Then a local optimal solution with an optimal value $1047294\$$ is found in 0.22 min. On the other hand, we solve this system by using MILP-IPM, where the parameters in step 1 are set the same as the case 2. After the optimization, a local optimal solution with an optimal value $1040676\$$ is found in 1.12 min. The results are compared with other methods in Table \ref{cost10withloss}.

As we can see in the Table \ref{cost10withloss}, in comparison with IPM, MILP-IPM can solve to a much better solution in a reasonable time. Meanwhile, the total generation cost obtained by MILP-IPM is lower than most of the results reported in the literatures.

\begin{table}[!ht]
\caption{\label{cost10withloss} Summary results for the 10-unit system with loss}
\centering
\resizebox{0.46\textwidth}{!}{
\begin{tabular}{ccccc}
\hline
\multirow{2}{*}{Method} & \multicolumn{3}{c}{Total generation cost (\$)}&\multirow{2}{*}{S-time(min)} \\
\cline{2-4}
& Minimum & Average & Maximum & \\
\hline
EP\cite{EPSQP2002}        &  1054685&  1057323&      NA&    47.23\\
EP-SQP\cite{EPSQP2002}    &  1052668&  1053771&      NA&    27.53\\
GA\cite{ABC2011}          &  1052251&  1058041& 1062511&     4.59\\
MHEP-SQP\cite{MHEPSQP2005}&  1050054&  1052349&      NA&    24.33\\
DGPSO  \cite{DGPSO2005}   &  1049167&  1051725&      NA&     5.99\\
PSO\cite{ABC2011}         &  1048410&  1052092& 1057170&     5.45\\
IPSO \cite{IPSO2009}      &  1046275&  1048145&      NA&     0.15\\
AIS \cite{AIS2011}        &  1045715&  1047050& 1048431&     30.96\\
CE \cite{ECE2011}         &  1044051&  1045159&      NA&     0.60\\
ECE \cite{ECE2011}        &  1043989&  1044470&      NA&     0.64\\
ABC \cite{ABC2011}        &  1043381&  1044963& 1046805&     4.55\\
CDBCO\cite{CDBCO2014}     &  1042900&  1044700&      NA&     1.66\\
HIGA \cite{HIGA2013}      &  1041088&  1042980& 1044927&     4.75\\
TVAC-IPSO\cite{TVACIPSO2012}&1041066&  1042118& 1043625&     3.16\\
ICA \cite{ICA2012}        &  1040758&  1041665& 1043175&     NA\\
EBSO \cite{EBSO2013}      &  1038915&  1039188& 1039272&     0.17\\
CSO\cite{CSO2015}         &  1038320&  1039374& 1042518&     1.39\\
CSAPSO \cite{CSAPSO2011}  &  1038251&  1039543&      NA&     0.83\\
EAPSO\cite{EAPSO2012}     &  1037898&  1038109& 1038238&     2.88\\
IPM      &   1047294       &    && 0.22\\
MILP-IPM &   1040676       &    && 1.12 \\
\hline
\end{tabular}}
\end{table}

\begin{table*}[!ht]
\caption{\label{output10withloss} Dispatch results (MW) for the 10-unit system}
\centering
\resizebox{0.98\textwidth}{!}{
\begin{tabular}{c|cccccccccccccc}
\hline
$t$  & unit $1$&  unit $2$&  unit $3$&  unit $4$&  unit $5$& unit $6$&  unit $7$&  unit $8$&  unit $9$&  unit $10$& Loss & $\Delta P_{t}$ & $\Delta P_{t}$\cite{CSO2015}& $\Delta P_{t}$\cite{EAPSO2012} \\
\hline
1& 150.0000& 142.2665& 199.1034& 60.0000& 122.8666& 122.4498 & 129.5904 & 47.0000 & 20.0000 & 55.0000 & 12.2767 & 0.0000& 0.0030 & 0.0000 \\
\hline
 2& 150.0000& 222.2665& 238.5251& 60.0000& 73.0000& 122.4498 & 129.5904 & 55.3121 & 20.0000 & 55.0000 & 16.1439 & 0.0000& 0.0011 & 0.0001 \\
\hline
 3& 150.0000& 302.2665& 312.1763& 60.0000& 73.0000& 122.4498 & 99.5904 & 85.3121 & 20.0000 & 55.0000 & 21.7952 & 0.0001& 0.0054 & 0.0254 \\
\hline
 4& 226.6242& 312.1330& 297.3995& 60.0000& 122.8665& 122.4498 & 129.5904 & 85.3121 & 20.0000 & 55.0000 & 25.3756 & 0.0001& 0.0001 & 0.0001 \\
\hline
 5& 226.6242& 392.1330& 297.8235& 60.0000& 122.8666& 122.4498 & 129.5904 & 85.3121 & 20.0000 & 55.0000 & 31.7997 & 0.0001& 0.0022 & 0.0001 \\
\hline
 6& 303.2484& 396.7994& 297.3995& 60.0000& 172.7331& 144.3175 & 129.5904 & 85.3121 & 20.0000 & 55.0000 & 36.4004 & 0.0000& 0.0015 & 0.0001 \\
\hline
 7& 379.8726& 396.7994& 306.5561& 70.4152& 172.7331& 122.4498 & 129.5904 & 85.3121 & 22.0571 & 55.0000 & 38.7859 & 0.0001& 0.0087 & 0.1612 \\
\hline
 8& 379.8726& 396.7994& 302.5307& 120.4152& 172.7331& 122.4498 & 129.5904 & 85.3121 & 52.0571 & 55.0000 & 40.7605 & 0.0001& 0.0014 & 0.0000 \\
\hline
 9& 456.4968& 396.7994& 297.3995& 148.8859& 222.5997& 122.4498 & 129.5904 & 90.0000 & 52.0571 & 55.0000 & 47.2787 & 0.0001& 0.0023 & 0.0000 \\
\hline
10& 456.4968& 396.7994& 337.6226& 191.2457& 222.5997& 160.0000 & 129.5904 & 120.0000 & 52.0571 & 55.0000 & 49.4117 & 0.0000& 0.0066 & 0.0000 \\
\hline
11& 460.4568& 396.7994& 340.0000& 241.2457& 243.0000& 160.0000 & 129.5904 & 120.0000 & 52.0571 & 55.0000 & 52.1495 & 0.0001& 0.0011 & 0.0000 \\
\hline
12& 456.4968& 460.0000& 340.0000& 241.2457& 237.0093& 160.0000 & 129.5904 & 120.0000 & 80.0000 & 55.0000 & 59.3423 & 0.0001& 0.0023 & 0.0001 \\
\hline
13& 456.4968& 396.7994& 337.6226& 191.2457& 222.5997& 160.0000 & 129.5904 & 120.0000 & 52.0571 & 55.0000 & 49.4117 & 0.0000& 0.0005 & 0.0002 \\
\hline
14& 456.4968& 396.7994& 297.3995& 148.8859& 222.5997& 122.4498 & 129.5904 & 90.0000 & 52.0571 & 55.0000 & 47.2787 & 0.0001& 0.0044 & 0.0000 \\
\hline
15& 379.8726& 396.7994& 314.9562& 110.0000& 172.7331& 122.4498 & 129.5904 & 85.3121 & 50.0000 & 55.0000 & 40.7137 & 0.0001& 0.0167 & 0.0002 \\
\hline
16& 303.2484& 396.7994& 292.4895& 60.0000& 122.8665& 122.4498 & 129.5904 & 85.3121 & 20.0000 & 55.0000 & 33.7562 & 0.0001& 0.0025 & 0.0000 \\
\hline
17& 226.6242& 396.7994& 324.8902& 60.0000& 122.8666& 122.4498 & 129.5904 & 55.3121 & 20.0000 & 55.0000 & 33.5327 & 0.0000& 0.0113 & 0.0001 \\
\hline
18& 303.2484& 396.7994& 306.5303& 70.4152& 172.7331& 122.4498 & 129.5904 & 85.3121 & 22.0571 & 55.0000 & 36.1359 & 0.0001& 0.0048 & 0.0000 \\
\hline
19& 379.8726& 396.7994& 302.5307& 120.4152& 172.7331& 122.4498 & 129.5904 & 85.3121 & 52.0571 & 55.0000 & 40.7605 & 0.0001& 0.0047 & 0.0002 \\
\hline
20& 456.4968& 460.0000& 340.0000& 135.3837& 222.5997& 160.0000 & 129.5904 & 115.3121 & 52.0571 & 55.0000 & 54.4398 & 0.0000& 0.0038 & 0.0000 \\
\hline
21& 456.4968& 389.5329& 302.4598& 120.4152& 222.5997& 160.0000 & 129.5904 & 85.3121 & 50.0000 & 55.0000 & 47.4070 & 0.0001& 0.0060 & 0.0000 \\
\hline
22& 379.8726& 309.5329& 297.3995& 87.6310& 172.7331& 122.4498 & 129.5904 & 85.3121 & 20.0000 & 55.0000 & 31.5216 & 0.0002& 0.0070 & 0.0000 \\
\hline
23& 303.2484& 229.5329& 224.1642& 60.0000& 122.8666& 122.4498 & 129.5904 & 85.3121 & 20.0000 & 55.0000 & 20.1645 & 0.0001& 0.0075 & 0.0001 \\
\hline
24& 226.6242& 222.2665& 205.8519& 60.0000& 73.0000& 122.4498 & 129.5904 & 85.3121 & 20.0000 & 55.0000 & 16.0950 & 0.0001& 0.0001 & 0.0002 \\
\hline
\end{tabular}}
\end{table*}

\begin{table}[!ht]
\caption{\label{analysis} Solution differences between MILP and MILP-IP}
\centering
\resizebox{0.4\textwidth}{!}{
\begin{tabular}{ccccc}
\hline
\multicolumn{5}{c}{case 3: 5-unit system}\\
\hline
0       & (0,3]    & (3,6]    & (6,9]   & (7,12]  \\
66.67\% & 12.50\%  &  8.33\%  &  7.50\% &  5.00\% \\
\hline
\multicolumn{5}{c}{case 4: 10-unit system}\\
\hline
0       & (0,10]  &  (10,20]   &  (20,30]  & (30,40] \\
70.83\% & 15.00\% &    5.42\%  &   5.00\%  &  3.75\%  \\
\hline
\end{tabular}}
\end{table}

In Table \ref{output10withloss}, we provide the generation dispatch results obtained by MILP-IPM for the 10-unit system. In this table, the violations of power balance constraints for CSO\cite{CSO2015} and EAPSO\cite{EAPSO2012} are also calculated due to their available solutions. For a fair comparison, we use the outputs shown in the corresponding literatures to compute the corresponding $\Delta P_{t}$. From the Table \ref{output10withloss} we notice that, although CSO\cite{CSO2015} and EAPSO\cite{EAPSO2012} can obtain lower total generation costs than MILP-IPM, but MILP-IPM outperforms them in terms of the feasibility of the solution for DED-VPE. In fact, when the original outputs are adopted for our case, all the violations are less than 8e-6. It means that, the feasibility of the solution obtained by MILP-IPM can be strictly satisfied.

\subsection{Results analysis}

From the simulation process in this section, we observe that, in the 5-unit system with transmission loss, $66.67\%$ of the solution points obtained by MILP-IPM are the same as those obtained by the individual MILP and differences of the rest range from 0.0960 to 11.72. In the 10-unit system with transmission loss, $70.83\%$ of the solution points obtained by MILP-IPM are the same as those obtained by the individual MILP and differences of the rest range from 0.0001 to 39.77. The details with respect to the solution differences between MILP and MILP-IP for both cases are given in Table \ref{analysis}.

To some extent, this phenomenon indicates that after solving MILP formulation in step 1, an initial point which is close to a good optimal solution for DED-VPE is yielded, and then in step 2, IPM starts its search from this initial point and tunes to a good optimal solution where all the constraints are satisfied, which is consistent with the analysis in the subsection 4.2.

\section{Conclusion}

In this paper, a deterministic MILP-IPM is proposed to solve the non-convex and non-differentiable DED-VPE. To avoid the intractable non-differentiable characteristic of DED-VPE, we derive a differentiable NLP formulation for DED-VPE. Although the NLP formulation can be directly solved by IPM, but the optimization will easily trap in a poor local optima. Therefore, MILP is integrated to generate a good initial point. And then, IPM can be used to exploit a better local optima for DED-VPE. Comparing with the heuristics which are inherently stochastic, MILP-IPM results are much more stable and the feasibility of the solution can be strictly guaranteed. So, MILP-IPM as a deterministic optimization technique is very promising to apply to the practical problems when VPE type factors are considered.

%\section*{Acknowledgment}

\ifCLASSOPTIONcaptionsoff
  \newpage
\fi

\end{document}